\documentclass[11pt]{article}
\usepackage{fullpage}

\usepackage{amsmath,amssymb,amsthm}

\usepackage{hyperref}
\hypersetup{
    colorlinks=true,
    linkcolor=blue,
    citecolor=red,
    urlcolor=blue,
    pdfborder={0 0 0}
}

\usepackage{color}
\usepackage{makeidx}

\usepackage{makeidx}

 \makeatletter
 \@addtoreset{equation}{section}
 \makeatother

 \newtheorem{ittheorem}{Theorem}
 \newtheorem*{theorem*}{Theorem}
 \newtheorem*{corollary*}{Corollary}
 \newtheorem*{lemma*}{Lemma}
 \newtheorem{itlemma}{Lemma}
 \newtheorem{itproposition}{Proposition}
 \newtheorem{itdefinition}{Definition}

 \newtheorem{itremark}{Remark}
 \newtheorem{itclaim}{Claim}
 \newtheorem{itcorollary}{\bf Corollary}

 \newenvironment{theorem}{\addtocounter{equation}{1}
 \begin{ittheorem}}{\end{ittheorem}}

\newcommand{\din}{\partial^{in}}
\newcommand{\dout}{\partial^{\, out}}

\newcommand{\doutie}{\partial^{\, out,ext}_{\infty}}
\newcommand{\dini}{\partial^{\, in}_\infty}
\newcommand{\douti}{\partial^{\, out}_\infty}
\newcommand{\doute}{\partial^{\, out,ext}}

\newcommand{\tatop}[2]{\genfrac{}{}{0pt}{1}{#1}{#2}}

\newcommand{\XX}{{\underline{X}}}
\newcommand{\xx}{{\underline{x}}}

\def\T{\mathcal{T}^*}
\def\S{\mathcal{S}}
\def\U{\mathcal{U}}
\def\M{\mathcal{M}}

\def\supp{\text{supp\,}}
\def\afm{\overline{|f|\,}^{\lower2pt\hbox{$\scriptstyle M$}}}
\def\fm{\overline{f\,}^{\lower2pt\hbox{$\scriptstyle M$}}}
\def\fmn{\overline{f_N\,}^{\lower2pt\hbox{$\scriptstyle M$}}}
\def\fmne{^{\displaystyle \eta}\overline{f_N\,}^{\lower2pt\hbox{$\scriptstyle M$}}}
\def\pfm{\overline{\psi f\,}^{\lower2pt\hbox{$\scriptstyle M$}}}
\def\phfm{\overline{\phi_X f_N\,}^{\lower2pt\hbox{$\scriptstyle M$}}}

\def\uro{\smash{{U}^{\!\!\!\!\raise5pt\hbox{$\scriptstyle o$}}}}

 \def \ba {\begin{array}}
 \def \ea {\end{array}}

 \def \Z {{\mathbb Z}}
 \def \Zd {{\mathbb Z^d}}
 \def \Rd {{\mathbb R^d}}
 
 \def \Znd {{\mathbb Z_n^d}}
 \def \R {{\mathbb R}}
 
 \def \N {{\mathbb N}}

 \def \cS {{\cal S}}
 \def \cE {{\cal E}}
 \def \cF {{\cal F}}
 \def \cA {{\cal A}}

 \def \cM {{\cal M}}
 
 \def \cU {{\cal U}}

 \def \La {{\Lambda}}

\newcommand{\Ld}{{\mathbb{L}}^d}
\newcommand{\Ldi}{\smash{{\mathbb{L}}^{d,\infty}}}
\newcommand{\Edi}{{\mathbb{E}}^{d,\infty}}

\newcommand{\Ed}{{\mathbb{E}}^d}

\newcommand{\dinf}{\, {\rm d}_{\infty}\, }

\def\cA{{\cal A}}   
\def\cE{{\cal E}} \def\cF{{\cal F}}  
   
\def\cM{{\cal M}}   
  \def\cS{{\cal S}} 

\def\cU{{\cal U}}   
 
\def\tY{{\widetilde Y}}

\def\R{{\mathbb R}}

\def\N{{\mathbb N}}

\def\eqref#1{(\ref{#1})}

\def\cU{\mathcal{U}}

\def\cS{\mathcal{S}}

\def\cM{\mathcal{M}}

\def\cF{\mathcal{F}}
\def\cE{\mathcal{E}}

\def\cA{\mathcal{A}}

\makeindex

\begin{document}

\title{A martingale minimax exponential inequality for Markov chains}

 \author{
\begin{tabular}{c c }
\text{Rapha\"el Cerf\footnote{
\noindent
DMA, Ecole Normale Sup\'erieure,
CNRS, PSL Research University, 75005 Paris.}
\footnote{
\noindent Laboratoire de Math\'ematiques d'Orsay, Universit\'e Paris-Sud, CNRS, Universit\'e
Paris--Saclay, 91405 Orsay.}}\\
\end{tabular}
}

\maketitle



\begin{abstract}
	We prove a new inequality controlling the large deviations of the empirical
	measure of a Markov chain. 
	This inequality 
	is based on 
	the martingale used by Donsker and
	Varadhan and the minimax theorem. 
	It holds for convex sets and it requires to take
	an infimum over the starting point.
In the case of a compact space, this inequality is a
partial improvement of
the large deviations estimates of Donsker and Varadhan.
In the case of a non compact space, we condition on the event that the process visits
$n$ times a compact subset of the space and we still obtain a control on the exponential
scale.

\end{abstract}



\def\X{X}
\def\Y{Y}


There is an abundant literature devoted to 
obtaining
exponential inequalities 
for Markov chains. A typical goal is to generalize to the Markovian setting
the nice inequalities existing for i.i.d. random variables, like the Bernstein,
Hoeffding or Chernoff inequalities.
Of course, some hypothesis are required on the Markov chains in order to derive
an exponential inequality.
Without pretending to make an extensive review, 
let us mention some 
important and interesting works on this subject.
Le{\'o}n and Perron \cite{LEPE} have derived
Hoeffding bounds for discrete 
reversible Markov chains which have an optimal exponential rate.
Cattiaux and Guillin \cite{CG} have obtained deviations bounds with the 
help of general functional inequalities.
Kargin \cite{KA} derives a Bernstein--Hoeffding inequality from a spectral gap condition. His result has been recently improved by Naor, Rao and Regev \cite{NSO}.
Chazottes and Redig \cite{CR} have introduced some assumptions on the moments
of the coupling time.
Gao, Guillin and Wu \cite{GGW} use a method of transportation--information inequality to establish a type of Bernstein inequality. 
Dedecker and Gou\"ezel \cite{DEGO} prove an Hoeffding inequality for geometrically
ergodic Markov chains.
Bertail and Cl\'emen\c{c}on \cite{BECI}, Adamczak and Bednorz \cite{AB}
rely on regeneration
techniques.
Moulos and Anantharam \cite{MOU}
have obtained
optimal Chernoff and Hoeffding bounds for finite state Markov chains, with a constant prefactor.                               


We shall follow here a different strategy, based on a refinement of the martingale technique employed
by 
Donsker and Varadhan \cite{DV1}. 
This note is essentially a little {c}omplement to 
the classical paper 
\cite{DV1}.
We prove a new exponential inequality
controlling the large deviations of the empirical measure of a Markov chain. 
	This inequality is based on 
	two ingredients: the martingale used by Donsker and
	Varadhan and the minimax theorem 
	(notice that the minimax theorem has also been 
	used
	in the derivation of Chernoff's type inequalities for independent random variables, see \cite{PANO,DZ}).
	This inequality holds for convex sets and it requires to take
	an infimum over the starting point.
In the case of a compact state space, this inequality is a
partial improvement of
the large deviations estimates of Donsker and Varadhan.
In the case of a non compact space, we condition on the event that the process visits
$n$ times a compact subset of the space and we still obtain a control on the exponential
scale.
\smallskip

\noindent
{\bf The model.}
Let $\X$ be a separable metric space.
Let $(X_n)_{n\in \N}$
be a stationary Markov chain whose state space is $\X$.
We denote by $P_x$
and $E_x$ the probability and expectation for 
the Markov chain $(X_n)_{n\in \N}$
starting from $X_0=x$.
For $u$ a continuous bounded function from $\X$ to $\R$, we set
$$\pi u(x)\,=\,E_x\big(u(X_1)\big)\,.$$
This notation is imported from \cite{DV1}, where $\pi$ stands for the transition kernel of the Markov chain.
We suppose that
%
%
%
for any real--valued function $u$ which is bounded and continuous on $X$,
the map 
$x\mapsto\pi u(x)$
is continuous on $\X$.
\smallskip

We define the local time $L_n$ at time $n\geq 1$ by setting, for $A$ a Borel subset of $\X$,
\begin{equation}
L_n(A)\,=\,
\sum_{k=0}^{n-1}1_{\{\,X_k\in A\,\}}\,.
\end{equation}
The local time $L_n$ is a random measure on $X$.
Let $\Y$ be a compact subset of $\X$. For $n\geq 1$, we define
$$\tau(\Y,n)\,=\,
\inf\,\big\{\,k\geq 1:L_{k}(\Y)=n\,\big\}\,.$$
In case the Markov chain visits $A$ strictly
less than $n$ times, the random time $\tau(\Y,n)$ is infinite.
We set, for $A$ a Borel subset of $\Y$,
$$L_n^\Y(A)\,=\,
\frac{1}{n}
L_{\tau(\Y,n)}(A)\,=\,
\frac{1}{n}
\sum_{k=0}^{\tau(\Y,n)-1}1_{\{\,X_k\in A\,\}}\,,
$$
where the sum runs over all the non--negative
integers in the case that
$\tau(\Y,n)=\infty$.
We denote by
$\cM(\Y)$ the set 
of the Borel probability measures on $\Y$, endowed with the topology of
weak convergence. The space $\Y$ being compact, the space
$\cM(\Y)$ is also compact (see chapter 1, section 5 in \cite{BIL}).
Whenever $\tau(\Y,n)$ is finite,
the random measure $L_{n}^{\Y}$ is a random element of $\cM(\Y)$.
We denote by $\mathcal{\cS^*}(Y)$ the 
collection of the real--valued continuous 
functions $u$ defined on $\X$ which satisfy
\begin{equation}
\label{bounded}
\exists\,c_1,c_2\in\R\quad
\forall y\in\X\qquad
0\,<\,c_1\,\leq\,u(y)\,\leq\,c_2\,,
\end{equation}
\begin{equation}
\label{subhar}
\forall y\in\X\setminus \Y\qquad{u(y)}
\,\geq\,
\pi u(y)
\,.
\end{equation}
The condition~\eqref{subhar} means that $u$ is 
$\pi$--superharmonic outside of $\Y$.
We define finally
$$
\forall\mu\in \M(\Y)\qquad
I(\mu)\,=\,-\inf_{u\in\S^*(Y)}\,\int_{\Y}\,\ln\Big(\frac{\pi u}{u}\Big)(x)\,\mu(dx)\,.$$
Let $\tY$ be a Borel set containing $Y$ and all 
the possible exit points from $Y$, i.e., such that
$Y\subset \tY$ and
$$\forall y\in Y\qquad P_y(X_1\in \tY)\,=\,1\,.$$
We state next our general exponential inequality.
\begin{theorem*}
\label{GD}
For any
closed convex subset $C$ of 
$\cM(\Y)$, we have
$$\forall n\geq 1\qquad
\inf_{x\in\tY}\,P_x\big(
L_n^\Y \in C
\big)\,\leq\,
\exp\Big(
-n\,\inf_{\mu\in C}\,I(\mu)
\Big)
\,.$$
\end{theorem*}
\noindent
An interesting feature of this theorem is that it does not require that the space $\X$ itself
is compact. Without further assumptions, it might very well happen that the process starting
from $x\in\Y$ does not return $n$ times to $\Y$. 
Yet on the event that 
$\{\,L_n^\Y \in C\,\}$, we have necessarily that
	$\tau(\Y,n)<\infty$
	(because $C\subset \cM(Y)$).
An unpleasant feature of the inequality is the infimum over the starting point.
It could well be that, for some set $Y$, we have to choose $\tY=X$, and that for some Markov chains,
the infimum over $X$ is $0$. 
However there are certainly many cases where the inequality provides an interesting
control. For instance, if the diameter of the jumps of the Markov chain are bounded, then 
the set $\tY$ can be taken to be an enlargement of $Y$, and if in addition the space is finite--dimensional,
then $\tY$ would also remain compact.
A variant of the inequality was obtained for the symmetric random walk on 
the lattice $\Z^d$ in \cite{BC} and it was a crucial tool to study
the random walk penalized by its range. In this situation, we can take for $Y$ any finite 
subset of $\Z^d$ and for $\tY$ the points of $\Z^d$ at distance at most $1$ from $Y$.
In the case that the space $\X$ is already compact, we can apply the theorem with
the choice $\Y=\X$. In this case, $\tau(\X,n)=n$, $nL_n^\Y=L_n$ and we obtain the following
corollary.
\begin{corollary*}
\label{expgd}
Suppose that the space $\X$ is compact.
Let $\cU_1$ be the 
collection of the continuous positive functions $u$ defined on $\X$.
We define
\begin{equation}
	\label{acti}
\forall\mu\in \M(\X)\qquad
I(\mu)\,=\,-\inf_{u\in\U_1}\,\int_{\X}\,\ln\Big(\frac{\pi u}{u}\Big)(x)\,\mu(dx)\,.\end{equation}
For any
closed convex subset $C$ of 
$\cM(\X)$, we have
\begin{equation}
	\label{expi}
	\forall n\geq 1\qquad
\inf_{x\in\X}\,P_x\Big(
	\frac{1}{n}L_n
	\in C
\Big)\,\leq\,
\exp\Big(
-n\,\inf_{\mu\in C}\,I(\mu)
\Big)
\,.\end{equation}
\end{corollary*}
\noindent
This corollary is a partial strengthening of the 
large deviations upper bound
of Donsker and Varadhan~(see \cite{DV1}, theorem~$1$), which we restate next:
for any closed subset $C$ of $\M(\X)$, we have
\begin{equation}
	\label{ldu}
	\forall {x\in\X}\qquad
\limsup_{n\to\infty}\,\frac{1}{n}\ln
P_x\Big(
	\frac{1}{n}L_n
	\in C
\Big)\,\leq\,
-\,\inf_{\mu\in C}\,I(\mu)
\,,\end{equation}
where
the function $I$ is given
by formula~\eqref{acti}.
This is the action functional of the large deviations principle
proved by Donsker and Varadhan under additional 
hypothesis~\cite{DV1}.
Their work contains also several results on $I$, as well
as alternative formulas for $I$.
In their work \cite{DV1}, Donsker and Varadhan make several 
hypothesis on the transition kernel $\pi$ of the Markov chain,
however these hypothesis are important for the proof of the large
deviations lower bound,
they are not needed for the proof of the large deviations upper bound.
The upper bound~\eqref{ldu} of Donsker and Varadhan is better 
than inequality~\eqref{expi}
in the sense that it holds for any starting
point $x$ and any closed subset $C$ of $\cM(\X)$, not necessarily convex.
The inequality~\eqref{expi} is better in the sense that it is a genuine exponential inequality,
which holds for any value $n\geq 1$.
Its weak side is that we have to take the infimum over all starting points
and that it holds only for convex sets.

Let us explain how we could derive the corollary starting from
the results of Donsker and Varadhan \cite{DV1}.
It is a classical fact that the quantity of the lefthand side 
of~\eqref{expi} is supermultiplicative,
i.e., that
$$
\forall m,n\geq 1\qquad
\inf_{x\in\X}\,P_x\Big(
	\frac{1}{m+n}L_{m+n}
	\in C
\Big)\,\geq\,
\inf_{x\in\X}\,P_x\Big(
	\frac{1}{n}L_n
	\in C
\Big)\,\times\,
\inf_{x\in\X}\,P_x\Big(
	\frac{1}{m}L_m
	\in C
\Big)\,,
$$
see \cite{DS}, chapter IV on uniform large deviations.
Equivalently, the sequence
$$
-\frac{1}{n}\ln
\Big( \inf_{x\in\X}\,P_x\Big(
	\frac{1}{n}L_n
	\in C
\Big)\Big)
$$
is subadditive. 
Unfortunately,
we cannot apply directly  the classical subadditive lemma, 
because the sequence might become infinite.
So we rely instead of the following variant.
\begin{lemma*}
\label{lc}
Let $f:\mathbb N\setminus\{0\}\to [0,\infty]$ be a subadditive map, i.e.,
$$\forall m,n\in\mathbb N\setminus\{0\}
\qquad f(m+n)\leq f(m)+f(n)\,.$$
We have
$$\liminf_{n\to\infty}\,
\frac{f(n)}{n}\,=\,\inf_{n\geq 1}\, 
\frac{f(n)}{n}\,.$$
\end{lemma*}
\noindent
The proof is very simple.
By the subadditive property, 
$$\forall d,n\in\mathbb{N}\setminus\{0\}\qquad
\frac{f(dn)}{dn}
\,\leq\,\frac{f(n)}{n}\,$$
and we conclude by sending first $d$ to $\infty$ and then taking
the infimum over $n$ in
$\mathbb{N}\setminus\{0\}$.
That's it!
Applying the lemma, 
we conclude that
\[
\liminf_{n\to\infty}\,-\frac{1}{n}\ln
\Big( \inf_{x\in\X}\,P_x\Big(
		\frac{1}{n}L_n
		\in C
	\Big)
\Big)\,=\,
\inf_{n\geq 1}\,-\frac{1}{n}\ln
\Big( \inf_{x\in\X}\,P_x\Big(
		\frac{1}{n}L_n
		\in C
	\Big)
\Big)\,,
\]
or equivalently
\begin{equation}
	\label{suba}
\limsup_{n\to\infty}\,\frac{1}{n}\ln
\Big( \inf_{x\in\X}\,P_x\Big(
	\frac{1}{n}L_n
	\in C
	\Big)
\Big)\,=\,
\sup_{n\geq 1}\,\frac{1}{n}\ln
\Big( \inf_{x\in\X}\,P_x\Big(
	\frac{1}{n}L_n
	\in C
	\Big)
\Big)\,.
\end{equation}
Rewriting~\eqref{suba} as an inequality and
combining it with the inequality~\eqref{ldu},
we obtain the inequality~\eqref{expi}.
It should be noted, though, that the proof of
the large deviations upper bound~\eqref{ldu} involves
a delicate exchange between an infimum and a supremum, that
is a sort of minimax theorem (see formula~(2.5) and thereafter
in \cite{DV1}). So the previous argument looks more tortuous
than the proof of the theorem, that we start next.

The crucial ingredient to derive the relevant probabilistic estimates
on the Markov chain is the family of martingales used by Donsker and Varadhan,
which we define thereafter.
Let $u$ be a function belonging to $\cS^*(Y)$.
For $n\geq 0$, we set
$$M_n\,=\,\Big(
\prod_{k=0}^{n-1}\frac{u(X_k)}{\pi u(X_k)}\Big)\,u(X_n)\,.$$
We claim that,
under $P_x$,
the process
$(M_n)_{n\in\N}$ is a martingale.
We denote by 
$E_x$ the expectation for the Markov chain $(X_n)_{n\in \N}$
starting from $X_0=x$.
Let us fix
$n\geq 1$.
We have
\begin{align*}
E_x\big(M_n\,\big|\,X_0,\dots,X_{n-1}\big)\,&=\,
\Big(\prod_{k=0}^{n-1}\frac{u(X_k)}{\pi u(X_k)}\Big)\,
E_x\big(u(X_n)\,|\,X_{n-1}\big)\cr
\,&=\,\Big(\prod_{k=0}^{n-1}\frac{u(X_k)}{\pi u(X_k)}\Big)\,
\pi u(X_{n-1})\,=\,M_{n-1}\,.
\end{align*}
The random time
$T={\tau(\Y,n)}-1$
is a stopping time for
$(M_n)_{n\in\N}$. Indeed, for $t\geq 1$, the event 
$\{\,
T
=t\,\}$ is measurable 
with respect to $L_{t+1}$, hence also to
	$X_1,\dots,X_{t}$.
	We apply next the optional stopping theorem:
\begin{equation}
	\label{opt}
	\forall t\geq 0\qquad
E_x\big(
M_{t\wedge 
	{T}
}\big)
\,=\,E_x(M_0)\,=\,u(x)\,.\end{equation}
We are interested in controlling $L_{\tau(Y,n)}=L_{T+1}$, so
we express $M_n$ with the help of the local time $L_{n+1}$ as follows:
\begin{align}
\label{locti}
M_n\,
&=\,\exp\Big(\sum_{0\leq k\leq n}\ln
\frac{u(X_k)}{\pi u(X_k)}\Big)\,\pi u(X_n)\cr
&=\,\exp\Big(\int_{y\in\X}\ln
\frac{u(y)}{\pi u(y)}\,dL_{n+1}(y)\Big)\,\pi u(X_n)
\,.
\end{align}
We
bound from below the lefthand member of~\eqref{opt}
and we use the identity~\eqref{locti} to get
(recall that $u$ is non--negative, and so is $M_n$)
\begin{align}
\label{belo}
\nonumber
E_x\big(
M_{t\wedge T}\big)
&\geq\,
E_x\Big(
M_{t\wedge T}
1_{\{T<\infty\}}
\Big)\\
&=\,
E_x\Bigg(
\,\exp\Big(\int_{y\in\X}\ln
\frac{u(y)}{\pi u(y)}\,dL_{t\wedge T+1}
(y)\Big)\,\pi u(X_{t\wedge T})
\,
1_{\{T<\infty\}}
\Bigg)
\,.
\end{align}
Recall that the function~$u$ belongs
to $\cS^*(Y)$.
The $\pi$--superharmonicity of $u$ on $X\setminus Y$ implies that
\begin{equation}
	\label{inte}
\int_{y\in\X}\ln
\frac{u(y)}{\pi u(y)}\,dL_{t\wedge T+1}(y)
\,\geq\,
\int_{y\in\Y}\ln
\frac{u(y)}{\pi u(y)}\,dL_{t\wedge T+1}(y)
\,.
\end{equation}
Reporting inequality~\eqref{inte} in 
inequalities~\eqref{belo} and~\eqref{opt}, we get
\begin{equation}
	\label{fgh}
	u(x)\,\geq\,
E_x\Bigg(
\,\exp\Big(
\int_{y\in\Y}\ln
\frac{u(y)}{\pi u(y)}\,dL_{t\wedge T+1}(y)
\Big)
\,
\pi u(X_{t\wedge T})
\,
1_{\{T<\infty\}}
\Bigg)
\,.
\end{equation}
On the event $\{T<\infty\}$, we have
$$\displaylines{
\forall y\in \Y\qquad
\lim_{t\to\infty}\,
L_{t\wedge T+1}(y)\,=\,
L_{T+1}(y)\,=\,
nL_n^\Y(y)\,,\cr
\lim_{t\to\infty}\,
X_{t\wedge T}\,=\,
X_{T}
\,.}$$
Applying
Fatou's lemma, we deduce from~\eqref{fgh} that
\begin{equation*}
u(x)
\,\geq\,
\,E_x\Bigg(
\,\exp\Big(
n\int_{y\in\Y}\ln
\frac{u(y)}{\pi u(y)}\,dL^\Y_{n}(y)
\Big)\,
\pi u(X_T) \, 
1_{\{T<\infty\}}
\Bigg)\,.
\end{equation*}
Yet $X_{T}$ belongs to $\Y$, thus 
\begin{equation}
	\pi u(X_T) \,=\,E_{X_T}\big(u(X_1)\big)
	\, \geq\, \inf_{y\in Y}\,
	E_{y}\big(u(X_1)\big)
	\, \geq\, \inf_{y\in \tY}\,u(y)
\end{equation}
and
\begin{equation}
u(x)
\,\geq\,
\Big(\inf_{y\in \tY}\,u(y)\Big)
\,E_x\Bigg(
\,\exp\Big(
n\int_{y\in\Y}\ln
\frac{u(y)}{\pi u(y)}\,dL^\Y_{n}(y)
\Big)\,
1_{\{T<\infty\}}
\Bigg)\,.
\end{equation}
Since $u$ belongs to $\S^*(Y)$, then the infimum on the right is positive.
Taking now the infimum over $x\in \tY$, we obtain
\begin{equation}
	\label{inbd}
\inf_{x\in \tY}
\,E_x\Bigg(
\,\exp\Big(
n\int_{y\in\Y}\ln
\frac{u(y)}{\pi u(y)}\,dL^\Y_{n}(y)
\Big)\,
1_{\{T<\infty\}}
\Bigg)\,\leq\,1\,.
\end{equation}
We proceed by bounding from below the lefthand member 
	of \eqref{inbd}
as follows:
for any $x\in  \tY$,
\begin{multline}
	\label{foui}
\,E_x\Bigg(
\,\exp\Big(
n\int_{y\in\Y}\ln
\frac{u(y)}{\pi u(y)}\,dL^\Y_{n}(y)
\Big)\,
1_{\{T<\infty\}}
\Bigg)\,\cr \geq
\,E_x\Bigg(
\,\exp\Big(
n\int_{y\in\Y}\ln
\frac{u(y)}{\pi u(y)}\,dL^\Y_{n}(y)
\Big)\,
1_{\{L_n^\Y\in C\}}
\Bigg)\\
 \geq\,
\,\exp\Bigg(
n\inf_{\mu\in C}\,\,
\int_{y\in\Y}\ln
\frac{u(y)}{\pi u(y)}\,d\mu(y)
\Bigg)\,
P_x\big(L_n^\Y\in C\big)\,.
\end{multline}
Recall that $C\subset\M(\Y)$, so that
$L_n^\Y\in C$ implies automatically that $T<\infty$.
Let us define
\begin{equation}
\label{phitc}
\phi_n(C)\,=\,
\inf_{x\in \tY}\,
P_x\big(L_n^\Y\in C\big)\,.
\end{equation}
The previous inequalities~\eqref{inbd} and~\eqref{foui}
yield that
\begin{equation}
\label{cvb}
\phi_n(C)\,\leq\,
\,\exp\Bigg(
-n\inf_{\mu\in C}\,\,
\int_{y\in\Y}\ln
\frac{u(y)}{\pi u(y)}\,d\mu(y)
\Bigg)\,.
\end{equation}
This inequality holds for any function $u$ in $\cS^*(Y)$.
In order to get a functional which is convex,
we perform a change of functions and we set
$\phi=\ln u$. We denote by $\T$ the image of $\cS^*(Y)$ under this change of functions, i.e.,
$$\T\,=\,\big\{\,\ln u:u\in\cS^*(Y)\,\big\}\,.$$
We rewrite inequality~\eqref{cvb} as follows: for any $\phi\in\T$, 
\begin{equation}
\label{cvc}
\phi_n(C)\,\leq\,
\,\exp\Bigg(
-n
\inf_{\mu\in
C
}\,\,
\int_{y\in\Y}
-\ln
\big(
e^{-\phi(y)}
\pi e^\phi(y)
\big)
\,d\mu(y)\Bigg)\,.
\end{equation}
We define a map $\Phi$ on $\T\times \cM(\Y)$
by
$$\Phi(\phi,\mu)\,=\,
\int_{y\in\Y}
-\ln
\big(
e^{-\phi(y)}
\pi e^\phi(y)
\big)
\,d\mu(y)\,.
$$
Optimizing the previous inequality~\eqref{cvc} 
over the function $\phi$ in~$\T$, we get
\begin{equation}
\label{hui}
\phi_n(C)\,\leq\,
\,\exp\Bigg(
	-n\,\,\sup_{\phi\in\T}\inf_{\mu\in
C
}\,\,
\Phi(\phi,\mu)
\Bigg)\,.
\end{equation}
The map $\Phi$ is linear in $\mu$. We shall next prove that it is convex in $\phi$.
In fact, the convexity in $\phi$ is a consequence of the convexity of the functions
$t\in\R\mapsto \exp(t)$ and
$t\in\R^+\mapsto -\ln(t)$. The delicate point is to check that the domain~$\T$
of definition of
$\Phi$ is convex. 
Let $\phi,\psi$ belong to $\T$ and let $\alpha,\beta\in]0,1[$ be such that $\alpha+\beta=1$.
There exist $u,v\in\cS^*(Y)$ such that 
$\phi=\ln u$, $\psi=\ln v$, whence
$$\alpha\phi+\beta\psi\,=\,\ln\big(u^\alpha v^\beta)\,,$$
and we have to check that
$u^\alpha v^\beta$ is in $\cS^*(Y)$. 
The continuity and condition~\eqref{bounded} are straightforward,
the only delicate point is 
the $\pi$--superharmonicity~\eqref{subhar}.
Let $y$ be a fixed point in $\X\setminus \Y$. We apply the H\"older inequality to the random variables
$u^\alpha(X_1)$, $v^\beta(X_1)$, with respect to the measure $P_y$
	and with the exponents $p=1/\alpha$, $q=1/\beta$.
We obtain
$$
	E_y\Big( u^\alpha(X_1) v^\beta(X_1)\Big)
	\,\leq\,
	E_y
	\Big(u(X_1)\Big)^\alpha
	E_y
	\Big(v(X_1)\Big)^\beta
	\,=\,
	\big(\pi u(y)\big)^\alpha
	\big(\pi v(y)\big)^\beta
\,.
$$
We finally use the fact that $u,v$ are $\pi$--superharmonic in $\X\setminus\Y$ to conclude that
\begin{equation*}
\forall y\in\X\setminus \Y\qquad
	 u^\alpha(y)
	 v^\beta(y)
\,\geq\,
\pi u^\alpha v^\beta (y)
\,.
\end{equation*}

\noindent
Moreover the set $C$ is compact and convex.
For any fixed $\phi\in\T$, the map
$\mu\in\cM(Y)\mapsto \Phi(\phi,\mu)$ is linear and
continuous for the weak convergence.
For any fixed $\mu\in\cM(Y)$, the map
$\phi\in\T\mapsto \Phi(\phi,\mu)$ is convex.
We are in position to apply
the famous minimax theorem (see \cite{FA}):
\begin{equation}
\label{supinf}
\sup_{\phi\in\T}\inf_{\mu\in
C
}\,\,
\Phi(\phi,\mu)
\,=\,
\inf_{\mu\in C }
\sup_{\phi\in\T}
\,\,
\Phi(\phi,\mu)
\,.
\end{equation}
Reporting in~\eqref{hui}, we conclude that
\begin{equation}
\label{huj}
\phi_n(C)\,\leq\,
\,\exp\Big(
	-n\,\,
	\inf_{\mu\in C }\,
	I(\mu)
	\,\,
\Big)\,,
\end{equation}
where $I(\mu)$ is defined by
\begin{equation}
	I(\mu)
	\,=\,
\sup_{\phi\in\T}\,
\Phi(\phi,\mu)\,.
\end{equation}
Coming back to the set of functions $\cS^*(Y)$, we have
\begin{equation}
	I(\mu)
	\,=\,
\sup_{u\in\cS^*(Y)}\,
\int_{y\in\Y}\ln
\frac{u(y)}{\pi u(y)}\,d\mu(y)
	\,=\,
-\inf_{u\in\cS^*(Y)}\,
\int_{y\in\Y}\ln
\Big(\frac{\pi u}{u}\Big)(y)\,d\mu(y)
\,,
\end{equation}
and so we get the inequality stated in the theorem.
\bibliographystyle{plain}
\bibliography{bobe}
 \thispagestyle{empty}

\end{document}